\documentclass[]{article}
\begin{document}
\newtheorem{proposition}{Proposition}[section]
\newtheorem{definition}{Definition}[section]
\newtheorem{lemma}{Lemma}[section]

\title{\bf Lie-like Algebras (Superalgeras)}
\author{Keqin Liu\\Department of Mathematics\\The University of British Columbia\\Vancouver, BC\\
Canada, V6T 1Z2}
\date{January 5, 2008}
\maketitle

\begin{abstract} We introduce the notion of a Lie-like algebra$^{\diamond}$ 
(superalgebra$^{\diamond}$) for $\diamond\in\{\,^{1-st},\, ^{2-nd},\, ^{3-rd}\,\}$.\end{abstract}

By bundling a family of algebras, we introduce six new generalizations of Lie algebras which are called Lie-like algebras or Lie-like superalgebras in this paper.  The examples of the new generalizations of Lie algebras come from an invariant algebra (or an invariant superalgebra), where the notion of an invariant algebra was introduced in Chapter 1 of \cite{Liu3}, and an invariant superalgebra is the counterpart of an invariant algebra in the context of superalgebras. Since an invariant algebra (or invariant superalgebra) can be regarded as an associative subalgebra of the associative algebra of linear transformations over a vector space
(or super vector space), the passage from an associative algebra to a Lie algebra can be extended to the passage from an invariant algebra (or an invariant superalgebra) to one of the new generalizations of Lie algebras. This paper consists of six sections. In the first three sections, we introduce Lie-like algebras$^{\diamond}$  with  $\diamond\in\{\,^{1-st},\, ^{2-nd},\, ^{3-rd}\,\}$. The counterparts of Lie-like algebras$^{\diamond}$ with  $\diamond\in\{\,^{1-st},\, ^{2-nd},\, ^{3-rd}\,\}$ in the context of superalgebras are called Lie-like superalgebras$^{\diamond}$ and introduced in the last three sections.

\medskip
Throughout this paper, we make the following conventions.
\begin{itemize}
\item All vector spaces are the vector spaces over a field $\bf{k}$.
\item Let $L$ be a vector space. We assume that a binary operation $L\times L\to L$ is always bilinear.
\item A vector space $V$ is called a super vector space if $V=V_0\oplus V_1$ is the direct sum of its two subspaces $V_0$ and $V_1$, where $0$, $1\in {\bf Z}_2={\bf Z}/2{\bf Z}$. The subspaces
$V_0$ and $V_1$ are called the even part and the odd part of $V$ respectively. 
\item A subspace $I$ of a super vector space $V=V_0\oplus V_1$ is called a  graded subspace of $V$ if $I=(I\cap V_0)\oplus (I\cap V_1)$.
\item Let $V=V_0\oplus V_1$ and $W=W_0\oplus W_1$ be two super vector spaces, and let $f: V\to W$ be a linear map. $f$ is called an even linear map if $f(V_{\alpha})\subseteq W_{\alpha}$ for all $\alpha\in {\bf Z}_2$, and $f$ is called an odd linear map if 
$f(V_{\alpha})\subseteq W_{\alpha +1}$ for all $\alpha\in {\bf Z}_2$.
\item If $V=V_0\oplus V_1$ is a super vector space, then $End\,V=End_0\, V\oplus End_1\, V$ is a super vector space, where the even part $End_0\, V$ and the odd part $End_1\, V$ are given by
$$
End_{0}\, V:=
\{\, f\in End\,V \,|\, \mbox{$f$ is an even linear map}\,\}
$$
and 
$$
End_{1}\, V:=
\{\, f\in End\,V \,|\, \mbox{$f$ is an odd linear map}\,\}.
$$
\end{itemize}

\medskip
\section{Lie-like Algebras$^{1-st}$}

Motivated by the properties of Hu-Liu square brackets in Proposition 1.4 of \cite{Liu3}, we  introduce the notion of a Lie-like algebra$^{1-st}$ in this section.

\begin{definition}\label{def1.1} Let $S$ be a nonempty set. A vector space $L$ is called a 
{\bf Lie-like algebra$^{1-st}$} induced by the set $S$ if there exists a family of binary operations
$$
 \Big\{\, [\; , \;]_k \;\Big |\;\mbox{$[\; , \;]_k : L\times L\to L$ is a binary map and $k\in S$} \,\Big\}
$$
such that the following two properties hold.
\begin{description}
\item[(i)] For $k\in S$, the binary operation $[\; , \;]_k$ is {\bf anti-symmetric}; that is
\begin{equation}\label{eq1}
 [x, y]_k=-[y, x]_k\quad\mbox{for $x$, $y\in L$ and $k\in S$}.
\end{equation}
\item[(ii)] For $x$, $y$, $z\in L$ and $h$, $k\in S$, the {\bf Jacobi-like identity$^{1-st}$} holds; that is,
\begin{equation}\label{eq2}
[[x, y]_h, z]_k+[[y, z]_k, x]_h+[[z, x]_h, y]_k=0.
\end{equation}
\end{description}
\end{definition}

A Lie-like algebra$^{1-st}$ $L$ induced by a set $S$ is also denoted by 
$(\, L\, , [\; , \;]_{k\in S}\,)$. Clearly, if $L$ is a Lie-like algebra$^{1-st}$ induced by a set $S$, then $L$ is a Lie algebra with respect to the binary operation $[\; , \;]_k$ for each $k\in S$. Hence, a Lie-like algebra$^{1-st}$ induced by a set $S$ is a bundle of Lie algebras satisfying the 
Jacobi-like identity$^{1-st}$. 

\medskip
A Lie algebra $(L, [\, ,\,])$ can be regarded as a Lie-like algebra$^{1-st}$ induced by a set 
$S$ if we define 
\begin{equation}\label{eq0}
 [x, y]_k=\phi (k)[x, y]\quad\mbox{for $x$, $y\in L$ and $k\in S$},
\end{equation}
where $\phi : k\mapsto \phi (k)$ is a map from $S$ to the field $\bf{k}$. We say that a Lie-like algebra$^{1-st}$ $(\, L\, , [\; , \;]_{k\in S}\,)$ is {\bf non-trivial} if (\ref{eq0}) does not hold for any binary operation $[\, , \,]: L\times L\to L$ and any map $\phi : S\to {\bf k}$. Our examples of Lie-like algebra$^{1-st}$ coming from invariant algebras are non-trivial.

\medskip
Let $L$ be a Lie-like algebra$^{1-st}$ $L$ induced by a set $S$. A subspace $I$ of  $L$ is called an {\bf ideal} of $L$ if $[L, \; I]_s\subseteq I$ for all $s\in S$. The following subspace of 
$L$ 
\begin{equation}\label{eq3}
L^{ann}_{1-st}:=\sum_{x, y\in L\atop h, k\in S}{\bf k}([x, y]_h-[x, y]_k)
\end{equation}
satisfies
\begin{equation}\label{eq4}
[L, \; L^{ann}_{1-st}]_s\subseteq L^{ann}_{1-st} \quad\mbox{for all $s\in S$},
\end{equation}
which proves that $L^{ann}_{1-st}$ is an ideal of $L$. The ideal $L^{ann}_{1-st}$ is called the
{\bf annihilator} of $L$, and the quotient space 
$L^{lie}_{1-st}:=\displaystyle\frac{L}{L^{ann}_{1-st}}$ is called the {\bf Lie-factor} of $L$. By 
(\ref{eq4}), the Lie factor $L^{lie}_{1-st}$ is a Lie algebra with respect to the bracket operation $[\; ,\;]$ defined by
\begin{equation}\label{eq5}
[\overline{x}, \;\overline{y}]: =\overline{[x, y]_k}, 
\end{equation}
where $x$, $y\in L$, $\overline{x}:=x+L^{ann}_{1-st}$ and $k\in S$.

\begin{definition}\label{def1.2} Let $i\in \{\, 2, \,3\,\}$. 
A Lie-like algebra$^{1-st}$ $L$ induced by a set $S$ is called a {\bf $i$-simple Lie-like algebra$^{1-st}$} if $L$ has no ideals which are not in the set 
$\{\, \{0\}, \,L^{ann}_{1-st},\, L\,\}$ and the set $\{\, \{0\}, \,L^{ann}_{1-st},\, L\,\}$
consists of $i$ distinct ideals.
\end{definition}

We now define a class of modules over a Lie-like algebra$^{1-st}$ $L$ induced by a set $S$.

\begin{definition}\label{def1.3} Let $L$ be a Lie-like algebra$^{1-st}$ induced by a set $S$. A vector space $V$ is called an {\bf ordinary module} over $L$ (or an {\bf ordinary $L$-module}) if there exist a family of linear maps
$$
 \Big\{\, f_k \;\Big |\;\mbox{$f_k : x\mapsto f_k(x)$ is a linear map from $L$ to $End(V)$ for $x\in L$ and $k\in S$} \,\Big\}
$$
such that
\begin{equation}\label{eq6}
 f_h([x, y]_k)=f_h(x)f_k(y)-f_h(y)f_k(x)
\end{equation}
and
\begin{equation}\label{eq7}
f_k(x)f_h(y)=f_h(x)f_k(y)
\end{equation}
for $x$, $y\in L$ and $h$, $k\in S$.
\end{definition}

An  ordinary module $V$ over a Lie-like algebra$^{1-st}$ induced by a set $S$ is also denoted by 
$(V, \{f_k\}_{k\in S})$. If $L$ is a Lie-like algebra$^{1-st}$ induced by a set $S$, then 
$(L, \{ad_k\}_{k\in S})$ is an ordinary module over $L$, where $ad_k$ is the {\bf adjoint map} defined by
\begin{equation}\label{eq8}
ad_k(x)(a):=[x, \; a]_k \quad\mbox{for $x$, $a\in L$ and $k\in S$}.
\end{equation}
$(L, \{ad_k\}_{k\in S})$ is called the {\bf adjoint module} over the Lie-like algebra$^{1-st}$ $L$ induced by the set $S$.

\medskip
Let $(V, \{f_k\}_{k\in S})$ be an ordinary module $V$ over a Lie-like algebra$^{1-st}$ $L$ induced by a set $S$. A subspace $U$ of $V$ is called an {\bf ordinary submodule} of $V$ if $f_k(x)(U)\subseteq U$ for 
$x\in L$ and $k\in S$. The subspace
\begin{equation}\label{eq9}
V^{ann}_{1-st}:=\sum_{x\in L,\; v\in V\atop h, k\in S}{\bf k}(f_h(x)-f_k(x))(v)
\end{equation}
is called the {\bf annihilator} of $V$. Clearly, the annihilator 
$V^{ann}_{1-st}$ of $V$ is an ordinary submodule of $V$.

\begin{definition}\label{def1.4} Let $i\in \{\, 2, \,3\,\}$. An ordinary module $V$ over a Lie-like algebra$^{1-st}$ $L$ induced by a set $S$ is said to be {\bf i-irreducible} if $V$ has no ordinary submodules which are not in the set 
$\{\, \{0\}, \,V^{ann}_{1-st},\, V\,\}$ and the set $\{\, \{0\}, \,V^{ann}_{1-st},\, V\,\}$
consists of $i$ distinct ordinary submodules.
\end{definition}

\medskip
We finish this section with the following remark.

\medskip\noindent
{\bf Remark} By Proposition 1.4 of \cite{Liu3}, there exist a bundle of Lie algebras
$$
\{\; (L, \; [\; , \;]_k) \;|\; \mbox{$L$ is a Lie algebra with the bracket $[\; , \;]_k$ for 
$k\in S$ }\;\}
$$
such that the following {\bf long Jacobi-like identity$^{1-st}$} holds
$$
[[x, y]_h, z]_k+[[y, z]_h, x]_k+[[z, x]_h, y]_k
+[[x, y]_k, z]_h+[[y, z]_k, x]_h+[[z, x]_k, y]_h=0
$$
for $x$, $y$, $z\in L$ and $h$, $k\in S$.

\medskip
\section{Lie-like Algebras$^{2-nd}$ }

Motivated by the properties of Hu-Liu angle brackets in Proposition 1.3 of \cite{Liu3}, we now introduce the notion of a Lie-like algebra$^{2-nd}$.

\begin{definition}\label{def2.1} Let $S$ be a nonempty set. A vector space $L$ is called a 
{\bf Lie-like algebra$^{2-nd}$} induced by the set $S$ if there exist a family of binary operations
$$
 \Big\{\, \langle\; , \;\rangle _k \;\Big |\;\mbox{$\langle\; , \;\rangle _k : L\times L\to L$ is a binary map and $k\in S$} \,\Big\}
$$
such that both the {\bf Jacobi-like identity$^{2-nd}$} 
\begin{equation}\label{eq10}
\langle \langle x, y\rangle _k, z\rangle _h=\langle x, \langle y, z\rangle _h\rangle _k+
\langle \langle x, z\rangle _h, y\rangle _k
\end{equation}
and the following identity
\begin{equation}\label{eq10'}
\langle \langle x, y\rangle _k, z\rangle _h=\langle \langle x, y\rangle _h, z\rangle _k
\end{equation}
hold for $x$, $y$, $z\in L$ and $h$, $k\in S$.
\end{definition}

If $h=k$, then (\ref{eq10}) becomes the Leibniz identity which was used to introduce Leibniz algebras in \cite{Loday}. Hence, if $L$ is a Lie-like algebra$^{2-nd}$ induced by a set $S$, then $L$ is a Leibniz algebra with respect to the binary operation $\langle\; , \;\rangle _k$ for each $k\in S$. Thus, a Lie-like algebra$^{2-nd}$ induced by a set $S$ is a bundle of Leibniz algebras satisfying the Jacobi-like identity$^{2-nd}$ and the identity (\ref{eq10'}). 

\medskip
Let $L$ be a Lie-like algebra$^{2-nd}$ induced by a set $S$. A subspace $I$ of  $L$ is called an {\bf ideal} of $L$ if $\langle I, \; L\rangle _s\subseteq I$ and 
$\langle L, \; I\rangle _s\subseteq I$  for all $s\in S$. The following subspace of $L$ 
\begin{equation}\label{eq11}
L^{ann, +}_{2-nd}:=\sum_{x, y\in L\atop h, k\in S}{\bf k}
(\langle x, y\rangle _h+\langle y, x\rangle _k)
\end{equation}
satisfies
\begin{equation}\label{eq12}
\langle L^{ann,+}_{2-nd}, \; L\rangle _s+\langle L, \; L^{ann,+}_{2-nd}\rangle _s
\subseteq L^{ann,+}_{2-nd} 
\quad\mbox{for all $s\in S$},
\end{equation}
which proves that $L^{ann,+}_{2-nd}$ is an ideal of $L$. The ideal $L^{ann,+}_{2-nd}$ is called the
{\bf plus annihilator} of $L$, and the quotient space 
$L^{lie,+}_{2-nd}:=\displaystyle\frac{L}{L^{ann,+}_{2-nd}}$ is called the {\bf Lie-factor} of $L$, which is a Lie algebra with respect to the bracket operation $[\; ,\;]$ defined by
\begin{equation}\label{eq13}
[\overline{x}, \;\overline{y}]: =\overline{\langle x, y\rangle_k}, 
\end{equation}
where $x$, $y\in L$, $\overline{x}:=x+L^{ann,+}_{2-nd}$ and $k\in S$. It is also easy to check that the following subspace of $L$ 
\begin{equation}\label{eq14}
L^{ann, -}_{2-nd}:=\sum_{x, y\in L\atop h, k\in S}{\bf k}
(\langle x, y\rangle _h-\langle x, y\rangle _k)
\end{equation}
satisfies
\begin{equation}\label{eq15}
\langle L^{ann,-}_{2-nd}, \; L\rangle _s+\langle L, \; L^{ann,-}_{2-nd}\rangle _s
\subseteq L^{ann,-}_{2-nd} 
\quad\mbox{for all $s\in S$}
\end{equation}
and
\begin{equation}\label{eq16}
L^{ann,-}_{2-nd}\subseteq L^{ann,+}_{2-nd},
\end{equation}
which prove that $L^{ann,-}_{2-nd}$ is an ideal of $L$ contained in the plus annihilator of $L$. The ideal $L^{ann,-}_{2-nd}$ is called the
{\bf minus annihilator} of $L$, and the quotient space 
$L^{lei,-}_{2-nd}:=\displaystyle\frac{L}{L^{ann,-}_{2-nd}}$ is called the {\bf Leibniz-factor} of $L$, which is a Leibniz algebra with respect to the angle operation $\langle \; ,\;\rangle$ defined by
\begin{equation}\label{eq17}
\langle \overline{\overline{x}}, \;\overline{\overline{y}}\rangle: =
\overline{\overline{\langle x, y\rangle _k}}, 
\end{equation}
where $x$, $y\in L$, $\overline{\overline{x}}:=x+L^{ann,-}_{2-nd}$ and $k\in S$. 

\begin{definition}\label{def2.2} Let $i\in \{\, 2, \,3, \,4\,\}$. 
A Lie-like algebra$^{2-nd}$ $L$ induced by a set $S$ is called a {\bf $i$-simple Lie-like algebra$^{2-nd}$} if $L$ has no ideals which are not in the set 
$\{\, \{0\}, \,L^{ann,-}_{1-st},\, L^{ann,+}_{1-st},\, L\,\}$ and the set 
$\{\, \{0\}, \,L^{ann,-}_{1-st},\, L^{ann,+}_{1-st},\, L\,\}$ 
consists of $i$ distinct ideals.
\end{definition}

We now define a class of modules over a Lie-like algebra$^{2-nd}$ $L$ induced by a set $S$.

\begin{definition}\label{def2.3} Let $L$ be a Lie-like algebra$^{2-nd}$ induced by a set $S$. A vector space $V$ is called an {\bf ordinary module} over $L$ (or an {\bf ordinary $L$-module}) if there exist a family of linear maps
$$
 \bigg\{\, f_k, \; g_k \;\bigg |\; \begin{array}{l}\mbox{$f_k : x\mapsto f_k(x)$ and 
$g_k : x\mapsto g_k(x)$ are linear}\\\mbox{maps from $L$ to $End(V)$ for $x\in L$ and $k\in S$}
\end{array} \,\bigg\}
$$
such that
\begin{equation}\label{eq18}
 f_h(\langle x, y\rangle _k)=f_h(x)f_k(y)-f_k(y)f_h(x),
\end{equation}
\begin{equation}\label{eq19}
 g_h(\langle x, y\rangle _k)=g_h(x)f_k(y)-f_k(y)g_h(x),
\end{equation}
\begin{equation}\label{eq20}
g_k(x)g_h(y)=g_h(x)f_k(y)=g_k(x)f_h(y),
\end{equation}
and
\begin{equation}\label{eq20a}
f_k(x)f_h(y)=f_h(x)f_k(y),\quad f_k(x)g_h(y)=f_h(x)g_k(y)
\end{equation}
for $x$, $y\in L$ and $h$, $k\in S$.
\end{definition}

The notion of a module over a Leibniz algebra, which was introduced in \cite{Loday3}, is a special case of the notion of an ordinary module over a Lie-like algebra$^{2-nd}$. 
An ordinary  module $V$ over a Lie-like algebra$^{2-nd}$ induced by a set $S$ is also denoted by 
$(V, \{f_k\}, \{g_k\}_{k\in S})$. If $L$ is a Lie-like algebra$^{2-nd}$ induced by a set $S$, then 
$(L, \{-r_k\}, \{\ell _k\}_{k\in S})$ is an ordinary module over $L$, where $r_k$ and $\ell _k$ are the {\bf right multiplication} and {\bf left multiplication}, respectively; that is, 
\begin{equation}\label{eq21}
r_k(x)(a):=\langle a, \; x\rangle _k \quad\mbox{and}\quad
\ell _k(x)(a):=\langle x, \; a\rangle _k\quad\mbox{for $x$, $a\in L$ and $k\in S$}.
\end{equation}
This ordinary $L$-module $(L, \{-r_k\}, \{\ell _k\}_{k\in S})$ is called the {\bf adjoint module} over the Lie-like algebra$^{2-nd}$ $L$ induced by the set $S$.

\medskip
Let $(V, \{f_k\}, \{g_k\}_{k\in S})$ be an ordinary module $V$ over a Lie-like algebra$^{2-nd}$ $L$ induced by a set $S$. A subspace $U$ of $V$ is called an {\bf ordinary submodule} of $V$ if 
$$
f_k(x)(U)\subseteq U\quad\mbox{and}\quad g_k(x)(U)\subseteq U\quad\mbox{for 
$x\in L$ and $k\in S$}. 
$$
The subspace
\begin{equation}\label{eq22}
V^{ann,+}_{2-nd}:=\sum_{x\in L,\; v\in V\atop h, k\in S}{\bf k}(g_h(x)-f_k(x))(v)
\end{equation}
is called the {\bf plus annihilator} of $V$ and the subspace
\begin{equation}\label{eq23}
V^{ann,-}_{2-nd}:=\sum_{x\in L,\; v\in V\atop h, k\in S}{\bf k}(f_h(x)-f_k(x))(v)
+\sum_{x\in L,\; v\in V\atop h, k\in S}{\bf k}(g_h(x)-g_k(x))(v)
\end{equation}
is called the {\bf minus annihilator} of $V$. Both plus annihilator $V^{ann,+}_{2-nd}$ and
minus annihilator $V^{ann,-}_{2-nd}$ are ordinary submodules of $V$. Moreover, we have
\begin{equation}\label{eq24}
V^{ann,-}_{2-nd}\subseteq V^{ann,+}_{2-nd}.
\end{equation}

\begin{definition}\label{def1.4} Let $i\in \{\, 2, \,3,\, 4\,\}$. An ordinary module $V$ over a Lie-like algebra$^{2-nd}$ $L$ induced by a set $S$ is said to be {\bf i-irreducible} if $V$ has no ordinary submodules which are not in the set 
$\{\, \{0\}, \,V^{ann,-}_{2-nd}, \,V^{ann,+}_{2-nd},\, V\,\}$ and the set 
$\{\, \{0\}, \,V^{ann,-}_{2-nd}, \,V^{ann,+}_{2-nd},\, V\,\}$ consists of $i$ distinct ordinary submodules.
\end{definition}

\medskip
\section{Lie-like Algebras$^{3-rd}$ }

The following identity
\begin{equation}\label{eq25}
[[x, y], \sigma (z)]+[[y, z], \sigma(x)]+[[z, x], \sigma(y)]=0
\end{equation}
was used in \cite{LiuPhd} and \cite{Liuphy} to generalize the Witt algebra, where 
$[\; ,\;]$ is an anti-symmetric binary operation and $\sigma$ is a linear map. Later, the same identity (\ref{eq25}) was  used in \cite{HLS} to define Hom-Lie algebras. Using invariant algebras introduced in \cite{Liu3}, we can produce several ways of constructing an algebra 
$(A, \, \circ)$ such that the algebra $(A, \, \circ)$  satisfies the following identity:
\begin{equation}\label{eq26}
(x\circ y)\circ \sigma (z)=\mathring{\sigma}(x)\circ (y\circ z)+(x\circ z)\circ \check{\sigma}(y),
\end{equation}
where $x$, $y$, $z\in A$, $\sigma$, $\mathring{\sigma}$, $\check{\sigma}\in End (A)$, and the binary operation $\circ$ is either anti-symmetric or not anti-symmetric. Clearly, (\ref{eq25}) is a special case of (\ref{eq26}).
By bundling a family of algebras which satisfy the identity (\ref{eq26}), we now introduce the notion of a Lie-like algebra$^{3-rd}$.

\begin{definition}\label{def3.1} A vector space $L$ is called a 
{\bf Lie-like algebra$^{3-rd}$} if there exist a nonempty set $S$, three subsets $\bf \Sigma$, 
$\mathring{\bf \Sigma}$ and $\check{\bf \Sigma}$ of $End(L)$ and a family of binary operations
$$
 \Big\{\, \circ_k \;\Big |\;\mbox{$\circ_k : L\times L\to L$ is a binary map and $k\in S$} \,\Big\}
$$
such that the following {\bf Jacobi-like identity$^{3-rd}$} holds
\begin{equation}\label{eq27}
(x\circ_h y)\circ_k \sigma (z)=\mathring{\sigma}(x)\circ_h (y\circ_k z)+(x\circ_k z)\circ_h \check{\sigma}(y),
\end{equation}
for $x$, $y$, $z\in L$, $h$, $k\in S$, $\sigma\in {\bf \Sigma}$, $\mathring{\sigma}\in\mathring{\bf \Sigma}$ and
$\check{\sigma}\in\check{\bf \Sigma}$.
\end{definition}

A Lie-like algebra$^{3-rd}$ $(\, L, \, S, \, {\bf \Sigma},\, \mathring{\bf \Sigma},\, \check{\bf \Sigma})$ is said to be {\bf anti-symmetric} if the binary operation $\circ_k$ is anti-symmetric for all $k\in S$.

\medskip
\section{Lie-like Superlgebras$^{1-st}$}

In the remaining part of this paper, $L=L_0\oplus L_1$ always denotes a super vector space with the even part $L_0$ and the odd part $L_1$, where $0$, $1\in {\bf Z}_2$. 

\medskip
The super counterpart of a Lie-like algebra$^{1-st}$ is introduced in the following

\begin{definition}\label{def4.1} Let $S$ be a nonempty set. A super vector space 
$L==L_0\oplus L_1$ is called a 
{\bf Lie-like superalgebra$^{1-st}$} induced by the set $S$ if there exist a family of binary operations
$$
 \Big\{\, [\; , \;]_k \;\Big |\;\mbox{$[\; , \;]_k : L\times L\to L$ is a binary map and $k\in S$} \,\Big\}
$$
such that the following three properties hold.
\begin{description}
\item[(i)]  For $k\in S$, the super vector space $L=L_0\oplus L_1$ is a superalgebra with respect to the binary operation $[\; , \;]_k$; that is,
\begin{equation}\label{eq28}
[L_{\alpha} , L_{\beta}]_k \subseteq L_{\alpha+\beta}
\quad\mbox{for $\alpha$ and $\beta\in {\bf Z}_2$}.
\end{equation}
\item[(ii)] For $k\in S$, the binary operation $[\; , \;]_k$ is {\bf super anti-symmetric}; that is
\begin{equation}\label{eq29}
 [x_{\alpha}, y_{\beta}]_k=-(-1)^{\alpha\beta}[y_{\beta}, x_{\alpha}]_k
\end{equation}
where $\alpha$, $\beta\in {\bf Z}_2$, $x_{\alpha}\in L_{\alpha}$ and $y_{\beta}\in L_{\beta}$.
\item[(iii)] The {\bf super Jacobi-like identity$^{1-st}$} holds; that is,
\begin{equation}\label{eq30}
[[x_{\alpha}, y_{\beta}]_h, z_{\gamma}]_k=[x_{\alpha}, [y_{\beta}, z_{\gamma}]_k]_h
+(-1)^{\beta\gamma}[[x_{\alpha}, z_{\gamma}]_h, y_{\beta}]_k,
\end{equation}
where $\alpha$, $\beta$, $\gamma\in {\bf Z}_2$, $x_{\alpha}\in L_{\alpha}$, 
$y_{\beta}\in L_{\beta}$ and $z_{\gamma}\in L_{\gamma}$.
\end{description}
\end{definition}

Clearly, if $L$ is a Lie-like superalgebra$^{1-st}$ induced by a set $S$, then $(L, [\; , \;]_k)$  is a Lie superalgebra in \cite{Kac} for each $k\in S$. Hence, a Lie-like superalgebra$^{1-st}$ induced by a set $S$ is a bundle of Lie superalgebras satisfying the 
super Jacobi-like identity$^{1-st}$. 

\medskip
If $L_0$ and $L_1$ are the even part and odd part of a Lie-like superalgebra$^{1-st}$ $L$ induced by a set $S$ respectively, then the even part $L_0$ is a Lie-like algebra$^{1-st}$ induced by the set $S$ and the odd part $L_1$ is an ordinary module over the Lie-like algebra$^{1-st}$ $L_0$.

\medskip
Let $L=L_0\oplus L_1$ be a Lie-like superalgebra$^{1-st}$ induced by a set $S$. A graded subspace $I$ of  $L$ is called a {\bf graded ideal} of $L$ if $[L, \; I]_s\subseteq I$ for all $s\in S$. 
The following graded subspace of $L$ 
\begin{equation}\label{eq31}
L^{sann}_{1-st}:=\sum_{x\in L_\alpha, y\in L_\beta\atop h, k\in S, 
\alpha, \beta\in {\bf Z}_2}{\bf k}([x, y]_h-[x, y]_k)
\end{equation}
satisfies
\begin{equation}\label{eq32}
[L, \; L^{sann}_{1-st}]_s\subseteq L^{sann}_{1-st} \quad\mbox{for all $s\in S$},
\end{equation}
which proves that $L^{sann}_{1-st}$ is an ideal of $L$. The ideal $L^{sann}_{1-st}$ is called the
{\bf super annihilator} of $L$, and the quotient space 
$L^{super}_{1-st}:=\displaystyle\frac{L}{L^{sann}_{1-st}}$ is called the {\bf super factor} of $L$. By 
(\ref{eq32}), the super factor $L^{super}_{1-st}$ is a Lie superalgebra with respect to the bracket operation $[\; ,\;]$ defined by
\begin{equation}\label{eq33}
[\overline{x}, \;\overline{y}]: =\overline{[x, y]_k}, 
\end{equation}
where $x\in L_\alpha$, $y\in L_\beta$ and $\overline{x}:=x+L^{super}_{1-st}$ and $k\in S$.

\begin{definition}\label{def1.2} Let $i\in \{\, 2, \,3\,\}$. 
A Lie-like superalgebra$^{1-st}$ $L$ induced by a set $S$ is called a {\bf $i$-simple Lie-like superalgebra$^{1-st}$} if $L$ has no ideals which are not in the set 
$\{\, \{0\}, \,L^{sann}_{1-st},\, L\,\}$ and the set $\{\, \{0\}, \,L^{sann}_{1-st},\, L\,\}$
consists of $i$ distinct graded ideals.
\end{definition}

We now define a class of modules over a Lie-like superalgebra$^{1-st}$ $L$ induced by a set $S$.

\begin{definition}\label{def4.3} Let $L=L_0\oplus L_1$ be a Lie-like superalgebra$^{1-st}$ induced by a set $S$. A super vector space $V=V_0\oplus V_1$ is called an {\bf ordinary module} over $L$ (or a {\bf $L$-module}) if there exist a family of linear maps
$$
 \bigg\{\, f_k\;\bigg |\; \begin{array}{l}\mbox{$f_k : x\mapsto f_k(x)$ is an even linear map}\\\mbox{from $L$ to $End(V)$ for $x\in L$ and $k\in S$}
\end{array} \,\bigg\}
$$
such that
\begin{equation}\label{eq37'}
 f_h([x_{\alpha}, y_{\beta}]_k)=f_h(x_{\alpha})f_k(y_{\beta})-
(-1)^{\alpha\beta}f_h(y_{\beta})f_k(x_{\alpha})
\end{equation}
and
\begin{equation}\label{eq38'}
f_k(x_{\alpha})f_h(y_{\beta})=f_h(x_{\alpha})f_k(y_{\beta})
\end{equation}
for $x_{\alpha}\in L_{\alpha}$, $y_{\beta}\in L_{\beta}$, $\alpha$, $\beta\in {\bf Z}_2$ and $h$, $k\in S$.
\end{definition}

An  ordinary module $V=V_0\oplus V_1$ over a Lie-like superalgebra$^{1-st}$ induced by a set $S$ is also denoted by 
$(V=V_0\oplus V_1, \{f_k\}_{k\in S})$. If $L=L_0\oplus L_1$ is a Lie-like superalgebra$^{1-st}$ induced by a set $S$, then 
$(L=L_0\oplus L_1, \{ad_k\}_{k\in S})$ is an ordinary module over $L$, where $ad_k$ is the {\bf adjoint map} defined by
\begin{equation}\label{eq39'}
ad_k(x_{\alpha})(a_{\beta}):=[x_{\alpha}, \; a_{\beta}]_k \quad\mbox{for $x_{\alpha}\in L_{\alpha}$, $a_{\beta}\in L_{\beta}$, $\alpha$, $\beta\in {\bf Z}_2$  and $k\in S$}.
\end{equation}
$(L=L_0\oplus L_1, \{ad_k\}_{k\in S})$ is called the {\bf adjoint module} over $L$ induced by the set $S$.

\medskip
Let $(V=V_0\oplus V_1, \{f_k\}_{k\in S})$ be an ordinary module over a Lie-like 
super-\linebreak
algebra$^{1-st}$ $L=L_0\oplus L_1$ induced by a set $S$. A graded subspace $U$ of $V$ is called an {\bf ordinary submodule} of $V$ if $f_k(x)(U)\subseteq U$ for 
$x\in L$ and $k\in S$. The subspace
\begin{equation}\label{eq40'}
V^{sann}_{1-st}:=\sum_{x_{\alpha}\in L_{\alpha},\; v_{\beta}\in V_{\beta}\atop \alpha, \beta\in {\bf Z}_2,\,h, k\in S}
{\bf k}(f_h(x_{\alpha})-f_k(x_{\alpha}))(v_{\beta})
\end{equation}
is called the {\bf super annihilator} of $V$ which is an ordinary submodule of $V$.

\begin{definition}\label{def4.4} Let $i\in \{\, 2, \,3\,\}$. An ordinary module $V=V_0\oplus V_1$ over a Lie-like algebra$^{1-st}$ $L=L_0\oplus L_1$ induced by a set $S$ is said to be {\bf i-irreducible} if $V$ has no ordinary submodules which are not in the set 
$\{\, \{0\}, \,V^{sann}_{1-st},\, V\,\}$ and the set $\{\, \{0\}, \,V^{sann}_{1-st},\, V\,\}$
consists of $i$ distinct ordinary submodules.
\end{definition}

\medskip
We finish this section with the following remark.

\medskip\noindent
{\bf Remark} There exist a bundle of Lie superalgebras
$$
 \bigg\{\,(L, \; [\; , \;]_k) \;\bigg |\; \begin{array}{l}\mbox{$L=L_0\oplus L_1$ is a Lie superalgebra}\\\mbox{with the bracket $[\; , \;]_k$ for $k\in S$}
\end{array} \,\bigg\}
$$
such that the following {\bf long super Jacobi-like identity$^{1-st}$} holds
\begin{eqnarray*}
(-1)^{\gamma\alpha}[[x_{\alpha}, y_{\beta}]_h, z_{\gamma}]_k+
(-1)^{\alpha\beta}[[y_{\beta}, z_{\gamma}]_h, x_{\alpha}]_k+
(-1)^{\beta\gamma}[[z_{\gamma}, x_{\alpha}]_h, y_{\beta}]_k+&&\\
+(-1)^{\gamma\alpha}[[x_{\alpha}, y_{\beta}]_k, z_{\gamma}]_h+
(-1)^{\alpha\beta}[[y_{\beta}, z_{\gamma}]_k, x_{\alpha}]_h+
(-1)^{\beta\gamma}[[z_{\gamma}, x_{\alpha}]_k, y_{\beta}]_h&=&0
\end{eqnarray*}
for $\alpha$, $\beta$, $\gamma\in {\bf Z}_2$, $x_{\alpha}\in L_{\alpha}$, $y_{\beta}\in L_{\beta}$, $z_{\gamma}\in L_{\gamma}$ and $h$, $k\in S$.

\medskip
\section{Lie-like Superalgebras$^{2-nd}$ }

The super counterpart of a Lie-like algebra$^{2-nd}$ is introduced in the following

\begin{definition}\label{def5.1} Let $S$ be a nonempty set. A super vector space 
$L=L_0\oplus L_1$ is called a 
{\bf Lie-like superalgebra$^{2-nd}$} induced by the set $S$ if there exist a family of binary operations
$$
 \Big\{\, \langle\; , \;\rangle_k \;\Big |\;\mbox{$\langle\; , \;\rangle_k : L\times L\to L$ is a binary map and $k\in S$} \,\Big\}
$$
such that the following three properties hold.
\begin{description}
\item[(i)]  For $k\in S$, the super vector space $L=L_0\oplus L_1$ is a superalgebra with respect to the binary operation $\langle\; , \;\rangle_k$; that is,
\begin{equation}\label{eq38}
\langle L_{\alpha} , L_{\beta}\rangle _k \subseteq L_{\alpha+\beta}
\quad\mbox{for $\alpha$ and $\beta\in {\bf Z}_2$}.
\end{equation}
\item[(ii)] Both the {\bf super Jacobi-like identity$^{1-st}$}
\begin{equation}\label{eq40}
\langle \langle x_{\alpha}, y_{\beta}\rangle _h, z_{\gamma}\rangle _k=
\langle x_{\alpha}, \langle y_{\beta}, z_{\gamma}\rangle _k\rangle _h
+(-1)^{\beta\gamma}\langle \langle x_{\alpha}, z_{\gamma}\rangle _k, y_{\beta}\rangle _h
\end{equation}
and the following identity
\begin{equation}\label{eq40a}
\langle \langle x_{\alpha}, y_{\beta}\rangle _h, z_{\gamma}\rangle _k
=\langle \langle x_{\alpha}, y_{\beta}\rangle _k, z_{\gamma}\rangle _h
\end{equation}
hold for $\alpha$, $\beta$, $\gamma\in {\bf Z}_2$, $x_{\alpha}\in L_{\alpha}$, 
$y_{\beta}\in L_{\beta}$ and $z_{\gamma}\in L_{\gamma}$.
\end{description}
\end{definition}

If $h=k$, then (\ref{eq40}) becomes the Leibniz superidentity in \cite{AAO}. Hence, if $L$ is a Lie-like superalgebra$^{2-nd}$ induced by a set $S$, then $L$ is a Leibniz superalgebra with respect to the binary operation $\langle\; , \;\rangle _k$ for each $k\in S$. Hence, a Lie-like superalgebra$^{2-nd}$ induced by a set $S$ is a bundle of Leibniz superalgebras satisfying the 
super Jacobi-like identity$^{2-nd}$ and the identity (\ref{eq40a}). 

\medskip
Let $L=L_0\oplus L_1$ be a Lie-like superalgebra$^{2-nd}$ induced by a set $S$. A graded subspace $I$ of  $L$ is called an {\bf ideal} of $L$ if $\langle I, \; L\rangle _s\subseteq I$ and 
$\langle L, \; I\rangle _s\subseteq I$  for all $s\in S$. The following subspace of $L$ 
\begin{equation}\label{eq41}
L^{sann, +}_{2-nd}:=\sum_{x\in L_\alpha, y\in L_\beta\atop \alpha, \beta\in {\bf Z}_2,\, h, k\in S}{\bf k}
(\langle x, y\rangle _h+(-1)^{\alpha\beta}\langle y, x\rangle _k)
\end{equation}
satisfies
\begin{equation}\label{eq42}
\langle L^{sann,+}_{2-nd}, \; L\rangle _s+\langle L, \; L^{sann,+}_{2-nd}\rangle _s
\subseteq L^{ann,+}_{2-nd} 
\quad\mbox{for all $s\in S$},
\end{equation}
which proves that $L^{sann,+}_{2-nd}$ is an ideal of $L$. The ideal $L^{sann,+}_{2-nd}$ is called the
{\bf plus super annihilator} of $L$, and the quotient space 
$L^{super,+}_{2-nd}:=\displaystyle\frac{L}{L^{sann,+}_{2-nd}}$ is called the {\bf plus super factor} of $L$, which is a Lie superalgebra with respect to the bracket operation $\langle\; ,\;\rangle$ defined by
\begin{equation}\label{eq43}
[\overline{x}, \;\overline{y}]: =\overline{\langle x, y\rangle _k}, 
\end{equation}
where $x$, $y\in L$, $\overline{x}:=x+L^{sann,+}_{2-nd}$ and $k\in S$. It is also easy to check that the following subspace of $L$ 
\begin{equation}\label{eq44}
L^{sann, -}_{2-nd}:=\sum_{x\in L_\alpha, y\in L_\beta\atop \alpha, \beta\in {\bf Z}_2,\, h, k\in S}{\bf k}
(\langle x, y\rangle _h-\langle x, y\rangle _k)
\end{equation}
satisfies
\begin{equation}\label{eq45}
\langle L^{sann,-}_{2-nd}, \; L\rangle _s+\langle L, \; L^{sann,-}_{2-nd}\rangle _s
\subseteq L^{sann,-}_{2-nd} 
\quad\mbox{for all $s\in S$}
\end{equation}
and
\begin{equation}\label{eq46}
L^{sann,-}_{2-nd}\subseteq L^{sann,+}_{2-nd},
\end{equation}
which prove that $L^{sann,-}_{2-nd}$ is an ideal of $L$ contained in the plus super annihilator of $L$. The ideal $L^{sann,-}_{2-nd}$ is called the
{\bf minus super annihilator} of $L$, and the quotient space 
$L^{super,-}_{2-nd}:=\displaystyle\frac{L}{L^{ann,-}_{2-nd}}$ is called the {\bf minus super factor} of $L$, which is a Leibniz superalgebra with respect to the angle operation $\langle \; ,\;\rangle$ defined by
\begin{equation}\label{eq47}
\langle \overline{\overline{x}}, \;\overline{\overline{y}}\rangle: =
\overline{\overline{\langle x, y\rangle _k}}, 
\end{equation}
where $x$, $y\in L$, $\overline{\overline{x}}:=x+L^{sann,-}_{2-nd}$ and $k\in S$. 

\begin{definition}\label{def5.2} A Lie-like superalgebra$^{2-nd}$ $L$ induced by a set $S$ is called a \linebreak{\bf $i$-simple Lie-like superalgebra$^{2-nd}$} if $L$ has no ideals which are not in the set 
$\{\, \{0\}, \,L^{sann,-}_{1-st},\, L^{sann,+}_{1-st},\, L\,\}$ and the set 
$\{\, \{0\}, \,L^{sann,-}_{1-st},\, L^{sann,+}_{1-st},\, L\,\}$ 
consists of $i$ distinct ideals, where $i\in \{\, 2, \,3, \,4\,\}$.
\end{definition}

We now define a class of modules over a Lie-like superalgebra$^{2-nd}$ $L$ induced by a set $S$.

\begin{definition}\label{def5.3} Let $L=L_0\oplus L_1$ be a Lie-like superalgebra$^{2-nd}$ induced by a set $S$. A vector space $V=V_0\oplus V_1$ is called an {\bf ordinary module} over $L$ (or an {\bf ordinary $L$-module}) if there exist a family of linear maps
$$
 \bigg\{\, f_k, \; g_k \;\bigg |\; \begin{array}{l}\mbox{$f_k : x\mapsto f_k(x)$ and 
$g_k : x\mapsto g_k(x)$ are even linear}\\\mbox{maps from $L$ to $End(V)$ for $x\in L$ and $k\in S$}
\end{array} \,\bigg\}
$$
such that
\begin{equation}\label{eq48}
 f_h(\langle x_{\alpha}, y_{\beta}\rangle _k)=(-1)^{\alpha\beta}f_h(x_{\alpha})f_k(y_{\beta})-
f_k(y_{\beta})f_h(x_{\alpha}),
\end{equation}
\begin{equation}\label{eq49}
 g_h(\langle x_{\alpha}, y_{\beta}\rangle _k)=(-1)^{\alpha\beta}g_h(x_{\alpha})f_k(y_{\beta})-f_k(y_{\beta})g_h(x_{\alpha}),
\end{equation}
\begin{equation}\label{eq50}
g_k(x_{\alpha})g_h(y_{\beta})=g_h(x_{\alpha})f_k(y_{\beta})=g_k(x_{\alpha})f_h(y_{\beta})
\end{equation}
and
\begin{equation}\label{eq50a}
f_k(x_{\alpha})f_h(y_{\beta})=f_h(x_{\alpha})f_k(y_{\beta}),
\quad f_k(x_{\alpha})g_h(y_{\beta})=f_h(x_{\alpha})g_k(y_{\beta})
\end{equation}
for  $x_{\alpha}\in L_{\alpha}$, $y_{\beta}\in L_{\beta}$, $\alpha$, $\beta\in {\bf Z}_2$ and $h$, $k\in S$.
\end{definition}

An ordinary  module $V=V_0\oplus V_1$ over a Lie-like superalgebra$^{2-nd}$ induced by a set $S$ is also denoted by 
$(V=V_0\oplus V_1, \{f_k\}, \{g_k\}_{k\in S})$. If $L$ is a Lie-like superalgebra$^{2-nd}$ induced by a set $S$, then $L$ itself can be made into an ordinary module over $L$,
which is called the {\bf adjoint module} over the Lie-like superalgebra$^{2-nd}$ $L$ induced by the set $S$.

\medskip
Let $(V=V_0\oplus V_1, \{f_k\}, \{g_k\}_{k\in S})$ be an ordinary module over a Lie-like superalgebra$^{2-nd}$ $L=L_0\oplus L_1$ induced by a set $S$. A graded subspace $U$ of $V$ is called an {\bf ordinary submodule} of $V$ if 
$$
f_k(x)(U)\subseteq U\quad\mbox{and}\quad g_k(x)(U)\subseteq U\quad\mbox{for 
$x\in L$ and $k\in S$}. 
$$
The subspace
\begin{equation}\label{eq51}
V^{sann,+}_{2-nd}:=\sum_{x_{\alpha}\in L_{\alpha},\; v_{\beta}\in V_{\beta}\atop \alpha, \beta\in {\bf Z}_2,\,h, k\in S}{\bf k}(g_h(x_{\alpha})-f_k(x_{\alpha}))(v_{\beta})
\end{equation}
is called the {\bf plus super annihilator} of $V$ and the subspace
\begin{eqnarray}\label{eq52}
V^{sann,-}_{2-nd}:&=&\sum_{x_{\alpha}\in L_{\alpha},\; v_{\beta}\in V_{\beta}\atop \alpha, \beta\in {\bf Z}_2,\,h, k\in S}{\bf k}(f_h(x_{\alpha})-f_k(x_{\alpha}))(v_{\beta})+\nonumber\\
&&\quad +\sum_{x_{\alpha}\in L_{\alpha},\; v_{\beta}\in V_{\beta}\atop \alpha, \beta\in {\bf Z}_2,\,h, k\in S}{\bf k}(g_h(x_{\alpha})-g_k(x_{\alpha}))(v_{\beta})
\end{eqnarray}
is called the {\bf minus super annihilator} of $V$. Both plus super annihilator 
$V^{sann,+}_{2-nd}$ and
minus super annihilator $V^{sann,-}_{2-nd}$ are ordinary submodules of $V$. Moreover, we have
\begin{equation}\label{eq53}
V^{sann,-}_{2-nd}\subseteq V^{sann,+}_{2-nd}.
\end{equation}

\begin{definition}\label{def5.4} An ordinary module 
$V=V_0\oplus V_1$ over a Lie-like superalgebra$^{2-nd}$ $L=L_0\oplus L_1$ induced by a set $S$ is said to be {\bf i-irreducible} if $V$ has no ordinary submodules which are not in the set 
$\{\, \{0\}, \,V^{sann,-}_{2-nd}, \,V^{sann,+}_{2-nd},\, V\,\}$ and the set 
$\{\, \{0\}, \,V^{sann,-}_{2-nd}, \,V^{sann,+}_{2-nd},\, V\,\}$ consists of $i$ distinct ordinary submodules, where $i\in \{\, 2, \,3,\, 4\,\}$.
\end{definition}

\medskip
\section{Lie-like Superlgebras$^{3-rd}$ }

The super counterpart of a Lie-like superalgebra$^{3-rd}$ is introduced in the following

\begin{definition}\label{def6.1} A super vector space $L=L_0\oplus L_1$ is called a 
{\bf Lie-like super-\linebreak algebra$^{3-rd}$} if there exist a nonempty set $S$, three subsets $\bf \Sigma$, 
$\mathring{\bf \Sigma}$ and $\check{\bf \Sigma}$ of $End_0\,L$ and a family of binary operations
$$
 \Big\{\, \circ_k \;\Big |\;\mbox{$\circ_k : L\times L\to L$ is a binary map and $k\in S$} \,\Big\}
$$
such that the following two properties hold.
\begin{description}
\item[(i)]  For $k\in S$, the super vector space $L=L_0\oplus L_1$ is a superalgebra with respect to the binary operation $\circ_k$; that is,
\begin{equation}\label{eq54}
L_{\alpha} \circ_k L_{\beta} \subseteq L_{\alpha+\beta}
\quad\mbox{for $\alpha$ and $\beta\in {\bf Z}_2$}.
\end{equation}
\item[(ii)] The {\bf super Jacobi-like identity$^{3-rd}$} holds; that is,
\begin{equation}\label{eq55}
(x_{\alpha}\circ_h y_{\beta})\circ_k \sigma (z_{\gamma})=
\mathring{\sigma}(x_{\alpha})\circ_h (y_{\beta}\circ_k z_{\gamma})
+(-1)^{\beta\gamma}(x_{\alpha}\circ_k z_{\gamma})\circ_h \check{\sigma}(y_{\beta}),
\end{equation}
where $\alpha$, $\beta$, $\gamma\in {\bf Z}_2$, $x_{\alpha}\in L_{\alpha}$, 
$y_{\beta}\in L_{\beta}$, $z_{\gamma}\in L_{\gamma}$, $h$, $k\in S$, $\sigma\in {\bf \Sigma}$, $\mathring{\sigma}\in\mathring{\bf \Sigma}$ and
$\check{\sigma}\in\check{\bf \Sigma}$.
\end{description}
\end{definition}

A Lie-like superalgebra$^{3-rd}$ $(\, L=L_0\oplus L_1, \, S, \, {\bf \Sigma},\, \mathring{\bf \Sigma},\, \check{\bf \Sigma})$ is said to be {\bf super anti-symmetric} if the binary operation $\circ_k$ is super anti-symmetric for all $k\in S$.

\bigskip
We finish this paper with two remarks.

\medskip
{\bf Remark 1.} The six generalizations of Lie algebras in this paper are obtained by  bundling a family of algebras which belong to the same class of algebras. Other wilder generalizations of Lie algebras can be introduced by  bundling a family of algebras which do not belong to the same class of algebras.

\medskip
{\bf Remark 2.} Using ordinary modules is just one way of studying the representations of 
Lie-like algebra$^{\diamond}$ (superalgebra$^{\diamond}$) for $\diamond\in\{\,^{1-st},\, ^{2-nd}\,\}$. Another way of studying the representations of these Lie-like algebra objects can be introduced by using the ideas in Chapter 7 of \cite{Liu3}.

\end{document}